# DISCONTINUOUS GALERKIN ISOGEOMETRIC ANALYSIS FOR THE BIHARMONIC EQUATION

Stephen Edward Moore

*Katholische Hochschulgemeinde der Diözese Linz,*
*Petrinumstrasse 12/8, A-4040, Linz, Austria*

**Abstract**

We present and analyze an interior penalty discontinuous Galerkin Isogeometric Analysis (dG-IgA) method for the biharmonic equation in computational domain in $\mathbb{R}^d$ with $d = 2, 3$. The computational domain consist of several non-overlapping sub-domains or patches. We construct B-Spline approximation spaces which are discontinuous across patch interfaces. We present *a priori* error estimate in a discrete norm and numerical experiments to confirm the theory.

## 1. Introduction

In this paper, we consider the fourth-order Dirichlet boundary value problem: find $u : \overline{\Omega} \to \mathbb{R}$ such that

$$\Delta^2 u = f \quad \text{in } \Omega, \qquad u = g_0, \quad \text{and} \quad \mathbf{n} \cdot \nabla u = g_1 \quad \text{on} \quad \partial\Omega, \qquad (1)$$

where $\mathbf{n}$ is the external unit normal vector to the boundary $\partial\Omega$, the bi-Laplacian operator $\Delta^2 := \Delta\Delta$ with $\Delta$ as the Laplace operator, $f$ is a given source function, $g_0, g_1$ are boundary data and $\Omega \subset \mathbb{R}^d, d = 2, 3$ is a bounded Lipschitz computational domain with the boundary $\partial\Omega$. We assume that the domain $\Omega$ is generated by Computer Aided Design (CAD) system and represented by a single or multiple patches which are images of the parameter domain $(0, 1)^d$ by spline or NURBS maps.

The model problem (1) is an example of a fourth-order elliptic problem occurring usually in various model of computational mechanics such as the Bernoulli-Euler beam and the Poisson-Kirchhoff thin plate theories [27, 14]. Several numerical solution techniques for the fourth-order problem have been studied including conforming and non-conforming finite element methods (FEM)

---

*Email address:* `moorekwesi@gmail.com` (Stephen Edward Moore)

and mixed finite element methods see, e.g. [8, 28]. The construction of conforming methods for such problems require finite element spaces of $H^2(\Omega)$. Such $H^2-$conforming methods are known to require continuously differentiable (i.e. $C^1-$) piece-wise polynomials on the elements. This is however known to be considerably difficult to construct practically. Examples of such conforming finite elements for such a problem are the Argyris element which uses polynomials of degree $p = 5$ for triangular elements, the reduced Hsieh-Clough-Tocher (rHCT) or Hsieh-Clough-Tocher (HCT) element also called macro-elements, which uses cubic polynomials for sub-partition triangular elements and the Bogner-Fox-Schmit element which uses bi-cubic functions for rectangular elements . For the non-conforming finite element, a typical example for solving such a model problem is the Morley element which uses piece-wise quadratic polynomials, see e.g. [8].

Alternatively, the fourth-order partial differential equation (PDE) can be solved by using the interior penalty discontinuous Galerkin finite element methods. The interior penalty methods dates back to [11] where Douglas and Dupont combined conforming continuous finite element with penalty terms which led to consistent schemes to derive *a priori* error estimates. In [12], the continuous Galerkin (cG), discontinuous Galerkin (dG) and stabilization techniques combined in solving the fourth-order elliptic problems and applied to thin plate bending theory problems in structural mechanics and to a strain gradient theory problem. The continuous/discontinuous Galerkin method has further been applied to the biharmonic problem on closed surfaces [19]. A continuous interior penalty $hp-$version of the interior penalty discontinuous Galerkin finite element method for fourth-order elliptic problems has also been studied, see, e.g., [22, 25]. Finally, we mention that a continuous interior penalty method for fourth order elliptic boundary value problems including Kirchhoff plates on polygonal domains have been analyzed in [5, 6, 7]

In most recent times, isogeometric analysis (IgA) has been proposed to approximate solutions of PDEs, see, e.g. [3]. The IgA uses the same class of basis functions for both representing the geometry of the domain and approximating the solution of the PDEs. Furthermore, the IgA has $(p-1)-$continuous differentiable basis i.e. $C^{(p-1)}$ with degree $p \geq 1$ functions which makes it an ideal scheme for the approximation of higher order PDEs including the biharmonic problem (1), see, e.g. [26].

In this paper, we will present *a priori* error estimate for multi-patch interior penalty discontinuous Galerkin isogeometric analysis (dG-IgA) for biharmonic equation on conforming patches with matching meshes. The dG-IgA or Nitsche coupling method has been introduced and analyzed for second order elliptic problems, see e.g., [16, 17, 18, 21, 23]. Following the monograph of Di Pietro and Ern [10], our analysis will require three main ingredients namely; discrete stability, consistency and boundedness of the discrete bilinear form. Using approximation estimates for $h-$refined IgA meshes from [3] and [26], we will then present *a priori* error estimate in an appropriate discrete norm.

The rest of the paper is organized as follows. In Section 2, we introduce function spaces, weak formulation and the isogeometric analysis framework.

Section 3 involves the derivation of the interior penalty discontinuous Galerkin scheme. Then, in Section 4, we introduce a discrete NURBS space $V_h$ and a discrete norm $\|\cdot\|_h$ and prove the coercivity of the bilinear form. The boundedness of the bilinear form is asserted in a product space $V_{h,*} \times V_h$, where we will need another discrete norm $\|\cdot\|_{h,*}$ defined on the vector space $V_{h,*}$. The error analysis of the dG-IgA scheme is presented in Section 5. In Section 6, we present and discuss numerical experiments to confirm our theoretical results. Finally, we draw some conclusions and discuss future works in Section 7.

## 2. Preliminaries

Let $\Omega$ be a bounded Lipschitz domain with boundary $\partial\Omega$. We introduce the Sobolev space $H^s(\Omega) := \{v \in L_2(\Omega) : D^\alpha v \in L_2(\Omega), \text{ for } 0 \leq |\alpha| \leq s\}$, where $L_2(\Omega)$ denote the space of square integrable functions and let $\alpha = (\alpha_1, \ldots, \alpha_d)$ be a multi-index with non-negative integers $\alpha_1, \ldots, \alpha_d$, $|\alpha| = \alpha_1 + \ldots + \alpha_d$, $D^\alpha := \partial^{|\alpha|}/\partial x^\alpha$, and associate with the sobolev space $H^s(\Omega)$ the norm $\|v\|_{H^s(\Omega)} = \left(\sum_{0 \leq |\alpha| \leq s} \|D^\alpha v\|_{L_2(\Omega)}^2\right)^{1/2}$ see, e.g. [1].

The variational formulation of the biharmonic problem (1) reads: find $u \in V_0$ such that

$$a(u, v) = \ell(v), \quad \forall v \in V_0, \tag{2}$$

where the bilinear and linear forms are given by

$$a(u, v) = \int_\Omega \Delta u \Delta v \, dx \quad \text{and} \quad \ell(v) = \int_\Omega fv \, dx, \tag{3}$$

and the hyperplane and test space given by $V_D := \{v \in H^2(\Omega) : v = g_0, \; \mathbf{n} \cdot \nabla v = g_1 \text{ on } \partial\Omega\}$ and $V_0 := \{v \in H^2(\Omega) : v = 0, \; \mathbf{n} \cdot \nabla v = 0 \text{ on } \partial\Omega\}$. The existence and uniqueness of the variational problem (2) follows the well-known Lax-Milgram lemma see e.g. [8].

*2.1. B-Spline and Isogeometric Analysis*

We refer the reader to [9] for detailed study on B-splines or NURBS based Galerkin methods. For the unit interval $\widehat{\Omega} = [0, 1]$, we define a vector $\Xi = \{0 = \xi_1, \ldots, \xi_{n+p+1} = 1\}$ with a non-decreasing sequence of real numbers in the parameter domain $\widehat{\Omega} = [0, 1]$ called a knot vector. Given $\Xi, p \geq 1$, and $n$ the number of basis functions, the univariate B-spline basis functions are defined by the Cox-de Boor recursion formula as follows:

$$\widehat{B}_{i,0}(\xi) = \begin{cases} 1 & \text{if } \xi_i \leq \xi < \xi_{i+1}, \\ 0 & \text{else,} \end{cases}$$

$$\widehat{B}_{i,p}(\xi) = \frac{\xi - \xi_i}{\xi_{i+p} - \xi_i} \widehat{B}_{i,p-1}(\xi) + \frac{\xi_{i+p+1} - \xi}{\xi_{i+p+1} - \xi_{i+1}} \widehat{B}_{i+1,p-1}(\xi), \tag{4}$$

where a division by zero is defined to be zero. We note that a basis function of degree $p$ is $(p-m)$ times continuously differentiable across a knot value with the multiplicity $m$. If all internal knots have the multiplicity $m = 1$, then B-splines of degree $p$ are globally $(p-1)$−continuously differentiable.

In general for $d$−dimensional problems, the B-spline basis functions are tensor products of the univariate B-spline basis functions. Let $\Xi_\alpha = \{\xi_{1,\alpha}, \ldots, \xi_{n_\alpha+p_\alpha+1,\alpha}\}$ be the knot vectors for every direction $\alpha = 1, \ldots, d$. Let $\mathbf{i} := (i_1, \ldots, i_d), \mathbf{p} := (p_1, \ldots, p_d)$ and the set $\overline{\mathcal{I}} = \{\mathbf{i} = (i_1, \ldots, i_d) : i_\alpha = 1, 2, \ldots, n_\alpha; \ \alpha = 1, 2, \ldots, d\}$ be multi-indicies. Then the tensor product B-spline basis functions are defined by

$$\widehat{B}_{\mathbf{i},\mathbf{p}}(\xi) := \prod_{\alpha=1}^{d} \widehat{B}_{i_\alpha, p_\alpha}(\xi_\alpha), \tag{5}$$

where $\xi = (\xi_1, \ldots, \xi_d) \in \widehat{\Omega} = (0,1)^d$. The univariate and multivariate B-spline basis functions are defined in the parametric domain by means of the corresponding B-spline basis functions $\{\widehat{B}_{\mathbf{i},\mathbf{p}}\}_{\mathbf{i} \in \overline{\mathcal{I}}}$.

The distinct values $\xi_i, i = 1, \ldots, n$ of the knot vectors $\Xi$ provides a partition of $(0,1)^d$ creating a mesh $\widehat{\mathcal{K}}_h$ in the parameter domain where $\widehat{K}$ is a mesh element. The computational domain is described by means of a geometrical mapping $\mathbf{\Phi}$ such that $\Omega = \mathbf{\Phi}(\widehat{\Omega})$ and

$$\mathbf{\Phi}(\xi) := \sum_{\mathbf{i} \in \overline{\mathcal{I}}} C_{\mathbf{i}} \widehat{B}_{\mathbf{i},\mathbf{p}}(\xi), \tag{6}$$

where $C_{\mathbf{i}}$ are the control points.

In many practical applications, the computational domain $\Omega$ is decomposed into $N$ non-overlapping domain $\Omega_i$ called subdomains or patches denoted as $\mathcal{T}_h := \{\Omega_i\}_{i=1}^{N}$ such that $\overline{\Omega} = \bigcup_{i=1}^{N} \overline{\Omega}_i$ and $\Omega_i \cap \Omega_j = \emptyset$ for $i \neq j$. Each patch is the image of an associated geometrical mapping $\mathbf{\Phi}_i$ such that $\mathbf{\Phi}_i(\widehat{\Omega}) = \Omega_i, i = 1, \ldots, N$, see Figure 1. Let $F_{ij} = \partial\Omega_i \cap \partial\Omega_j, i \neq j$, denote the interior facets of two patches. The collection of all such interior facets is denoted by $\mathcal{F}_I$, and the collection of all Dirichlet facets $F_i = \partial\Omega_i \cap \partial\Omega$ is also denoted by $\mathcal{F}_D$. Furthermore, the collection of all internal and Dirichlet facets is denoted by $\mathcal{F} := \mathcal{F}_I \cup \mathcal{F}_D$.

We define the basis functions in the computational domain by means of the geometrical mapping as $B_{\mathbf{i},\mathbf{p}} := \widehat{B}_{\mathbf{i},\mathbf{p}} \circ \mathbf{\Phi}^{-1}$ and the discrete function space by

$$\mathbb{V}_h = \text{span}\{B_{\mathbf{i},\mathbf{p}} : \mathbf{i} \in \overline{\mathcal{I}}\}. \tag{7}$$

We assume that for each patch $\Omega_i, i = 1, \ldots, N$, the underlying mesh $\mathcal{K}_{h,i}$ is quasi-uniform i.e.

$$h_K \leq h_i \leq C_u h_K, \quad \text{for all} \quad K \in \mathcal{K}_{h,i}, \quad i = 1, \ldots, N, \tag{8}$$

where $C_u \geq 1$ and $h_i = \max\{h_K, K \in \mathcal{K}_{h,i}\}$ is the mesh size of $\Omega_i$ and $h_K$ is the diameter of of the mesh element $K$.

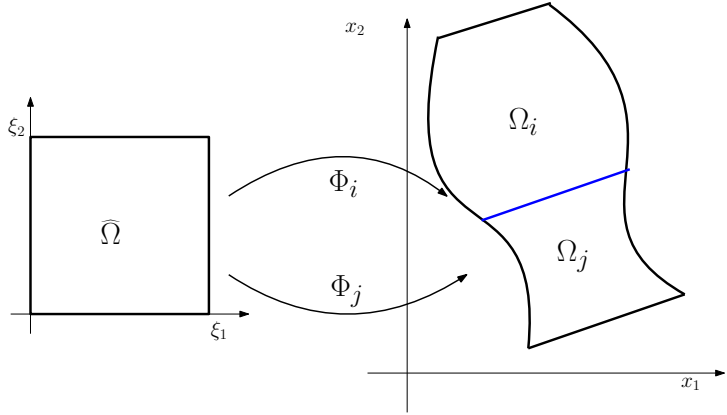

Figure 1: Illustration of the multi-patch isogeometric mapping.

## 3. Interior Penalty Variational Formulation

Let us recall some function spaces required for the derivation of interior penalty Galerkin schemes. We assign to each patch $\Omega_i$ an integer $s_i$ and collect them in the vector $\mathbf{s} = \{s_1, \ldots, s_N\}$. Let us now define the broken Sobolev space

$$H^{\mathbf{s}}(\Omega, \mathcal{T}_h) := \{v \in L_2(\Omega) : v|_{\Omega_i} \in H^{s_i}(\Omega_i), \ i = 1, \ldots, N\}, \tag{9}$$

and the corresponding broken Sobolev norm and semi-norm

$$\|v\|_{H^{\mathbf{s}}(\Omega, \mathcal{T}_h)} := \left(\sum_{i=1}^{N} \|v\|_{H^{s_i}(\Omega_i)}^2\right)^{1/2} \quad \text{and} \quad |v|_{H^{\mathbf{s}}(\Omega, \mathcal{T}_h)} := \left(\sum_{i=1}^{N} |v|_{H^{s_i}(\Omega_i)}^2\right)^{1/2}, \tag{10}$$

respectively.

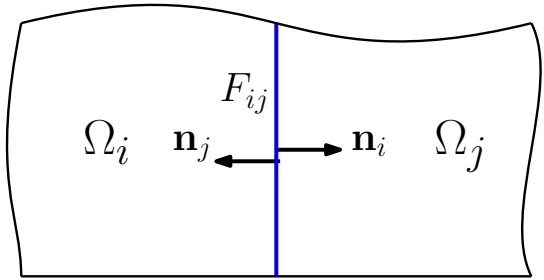

Figure 2: An illustration of the interior facet with normal vector $\mathbf{n}_i$ and $\mathbf{n}_j$.

Also, let $v_i$ and $v_j$ denote the restrictions of $v$ on patches $\Omega_i$ and $\Omega_j$, respectively. For the interior facets $F_{ij} \in \mathcal{F}_I$, let $\mathbf{n}_i$ be the outward unit normal vector with respect to $\Omega_i$ and let $\mathbf{n}_j$ be the outward unit normal vector with respect to $\Omega_j$, see Figure 2. We define the jump and average of the normal derivatives across the interior facets $F_{ij} \in \mathcal{F}_I$ of a smooth function $v \in H^2(\Omega, \mathcal{T}_h)$ by

$$[\![\nabla v]\!] := \mathbf{n}_i \cdot \nabla v_i + \mathbf{n}_j \cdot \nabla v_j, \quad and \quad \{\nabla v\} := \frac{1}{2}\left(\mathbf{n}_i \cdot \nabla v_i + \mathbf{n}_j \cdot \nabla v_j\right), \quad (11)$$

whereas the jump and average functions of the normal derivatives on the Dirichlet facets $F_i \in \mathcal{F}_D$ are given by

$$[\![\nabla v]\!] := \mathbf{n}_i \cdot \nabla v_i, \quad and \quad \{\nabla v\} := \mathbf{n}_i \cdot \nabla v_i. \quad (12)$$

Using the definitions of jumps and averages, the following equality

$$[\![ab]\!] = \{a\}[\![b]\!] + \{b\}[\![a]\!], \quad (13)$$

holds for some smooth functions $a : \Omega \to \mathbb{R}^2, b : \Omega \to \mathbb{R}$ on the interior facets $F_{ij} \in \mathcal{F}_I$.

The derivation of the interior penalty Galerkin scheme is given as follows; By integrating by parts over patch $\Omega_i, i = 1, \ldots, N$ for any $u \in H^4(\Omega_i)$ and $v \in H^2(\Omega_i)$, we obtain

$$\int_{\Omega_i} \Delta u \Delta v \, dx = \int_{\Omega_i} \Delta^2 uv \, dx - \int_{\partial \Omega_i} \Delta u (\mathbf{n} \cdot \nabla v) \, ds + \int_{\partial \Omega_i} \mathbf{n} \cdot \nabla \Delta u v \, ds. \quad (14)$$

Summing (14) over all patches $\Omega_i, i = 1, \ldots, N$, we obtain

$$\sum_{i=1}^N \int_{\Omega_i} \Delta u \Delta v \, dx = \sum_{i=1}^N \int_{\Omega_i} \Delta^2 uv \, dx - \sum_{i=1}^N \int_{\partial \Omega_i} \Delta u (\mathbf{n} \cdot \nabla v) \, ds$$
$$+ \sum_{i=1}^N \int_{\partial \Omega_i} \mathbf{n} \cdot \nabla \Delta u v \, ds. \quad (15)$$

Rewriting the boundary terms yields

$$\sum_{i=1}^N \int_{\partial \Omega_i} \Delta u (\mathbf{n} \cdot \nabla v) \, ds = \sum_{F_{ij} \in \mathcal{F}_I} \int_{F_{ij}} [\![\Delta u \nabla v]\!] \, ds + \sum_{F_i \in \mathcal{F}_D} \int_{F_i} \Delta u (\mathbf{n}_i \cdot \nabla v) \, ds, \quad (16)$$

and

$$\sum_{i=1}^N \int_{\partial \Omega_i} \mathbf{n} \cdot \nabla \Delta u v \, ds = \sum_{F_{ij} \in \mathcal{F}_I} \int_{F_{ij}} [\![\nabla \Delta u v]\!] \, ds + \sum_{F_i \in \mathcal{F}_D} \int_{F_i} \mathbf{n}_i \cdot \nabla \Delta u v \, ds. \quad (17)$$

Using the relation (13), we rewrite the terms on the interior facets $F_{ij}$ in (16) and (17) and substituted into (15) to obtain

$$\sum_{i=1}^N \int_{\Omega_i} fv \, dx = \sum_{i=1}^N \int_{\Omega_i} \Delta u \Delta v \, dx - \sum_{F_{ij} \in \mathcal{F}_I} \int_{F_{ij}} \{\nabla \Delta u\} [\![v]\!] \, ds$$

$$+ \sum_{F_{ij} \in \mathcal{F}_I} \int_{F_{ij}} \{\Delta u\} [\![\nabla v]\!] \, ds + \sum_{F_i \in \mathcal{F}_D} \int_{F_i} \mathbf{n}_i \cdot \nabla \Delta u v \, ds$$
$$+ \sum_{F_i \in \mathcal{F}_D} \int_{F_i} \Delta u (\mathbf{n}_i \cdot \nabla v) \, ds, \tag{18}$$

where we observe that for the exact solution $u$, we have $[\![\nabla \Delta u]\!] = 0$ and $[\![\Delta u]\!] = 0$ on the interior facets. Also $[\![u]\!] = 0$ and $[\![\nabla u]\!] = 0$ on $F \in \mathcal{F}$ whereas $[\![u - g_0]\!] = 0$ and $[\![\nabla u - g_1]\!] = 0$ on the Dirichlet facets $F_i \in \mathcal{F}_D$. We therefore add the following terms

$$- \sum_{F \in \mathcal{F}} \int_F \{\nabla \Delta v\} [\![u]\!] \, ds + \sum_{F \in \mathcal{F}} \int_F \{\Delta v\} [\![\nabla u]\!] \, ds.$$

Finally, we add a further consistency term that is zero and penalizes the interior facets by adding a penalty term to the bilinear form

$$\sum_{F \in \mathcal{F}} \int_F \frac{\delta_0}{h_i^3} [\![u]\!] [\![v]\!] \, ds + \sum_{F \in \mathcal{F}} \int_F \frac{\delta_1}{h_i} [\![\nabla u]\!] [\![\nabla v]\!] \, ds,$$

where the penalty parameters $\delta_0$ and $\delta_1$ are non-zero positive real numbers to be determined later.

The interior penalty variational scheme reads: find $u \in H^4(\Omega, \mathcal{T}_h)$ such that,

$$a_h(u, v) = \ell(v), \quad \forall v \in H^4(\Omega, \mathcal{T}_h), \tag{19}$$

where the bilinear form is given by

$$a_h(u, v) = \sum_{i=1}^{N} \int_{\Omega_i} \Delta u \Delta v \, dx$$
$$- \sum_{F \in \mathcal{F}} \int_F [\![\nabla v]\!] \{\Delta u\} \, ds + \sum_{F \in \mathcal{F}} \int_F [\![\nabla u]\!] \{\Delta v\} \, ds$$
$$+ \sum_{F \in \mathcal{F}} \int_F \{\nabla \Delta u\} [\![v]\!] \, ds - \sum_{F \in \mathcal{F}} \int_F \{\nabla \Delta v\} [\![u]\!] \, ds$$
$$+ \sum_{F \in \mathcal{F}} \frac{\delta_0}{h_i^3} \int_F [\![u]\!] [\![v]\!] \, ds + \sum_{F \in \mathcal{F}} \frac{\delta_1}{h_i} \int_F [\![\nabla u]\!] [\![\nabla v]\!] \, ds, \tag{20}$$

and the linear form reads as

$$\ell(v) = \sum_{i=1}^{N} \int_{\Omega_i} f v \, dx + \sum_{F_i \in \mathcal{F}_D} (\delta_0 v - \mathbf{n}_i \cdot \nabla \Delta v) g_0 \, ds$$
$$+ \sum_{F_i \in \mathcal{F}_D} (\delta_1 \mathbf{n}_i \cdot \nabla v + \Delta v) g_1 \, ds. \tag{21}$$

**Remark 1.** The penalty parameters $\delta_0$ and $\delta_1$ are known to depend on the B-spline or NURBS degree $p$ and the dimension $d$ of the domain $\Omega$, i.e. $\Omega \subset \mathbb{R}^d, d = 1, 2, 3$. For example, see [24] where the Dirichlet boundary value problem for Poisson's equation has been studied and an explicit expression of the penalty parameter obtained as $\delta_0 = \delta_1 = (p + 1)(p + d)/d$.

To show the consistency of the interior penalty Galerkin scheme, we require the weak continuity of the exact solution and its fluxes on the interior facets as follows;

**Lemma 1.** *The exact solution $u$ satisfies the following*

$$[\![u]\!] = 0, [\![\nabla u]\!] = 0, [\![\Delta u]\!] = 0, [\![\nabla \Delta u]\!] = 0, \quad \forall F \in \mathcal{F}.$$

PROOF. See [10, Lemma 4.3] and [25, Lemma 3].

**Theorem 1.** *Let $u \in H^4(\Omega)$ be the solution of (2). Then $u$ is the solution to dG-IgA variational identity (19).*

PROOF. By applying the integration by parts formula (14) and the relation (13), the variational form (19) for some $u \in H^4(\Omega, \mathcal{T}_h)$ yields

$$\begin{aligned}
0 = & \sum_{i=1}^{N} \int_{\Omega_i} \left(\Delta^2 u - f\right) v\, dx + \sum_{F \in \mathcal{F}} \int_F \left([\![\Delta u]\!]\{\nabla v\} + \{\Delta u\}[\![\nabla v]\!]\right) ds \\
& - \sum_{F \in \mathcal{F}} \int_F \left([\![\nabla \Delta u]\!]\{v\} + \{\nabla \Delta u\}[\![v]\!]\right) ds - \sum_{F \in \mathcal{F}} \int_F [\![\nabla v]\!]\{\Delta u\}\, ds \\
& + \sum_{F \in \mathcal{F}} \int_F [\![\nabla u]\!]\{\Delta v\}\, ds + \sum_{F \in \mathcal{F}} \int_F \{\nabla \Delta u\}[\![v]\!]\, ds - \sum_{F \in \mathcal{F}} \int_F \{\nabla \Delta v\}[\![u]\!]\, ds \\
& + \sum_{F \in \mathcal{F}} \frac{\delta_0}{h_i^3} \int_F [\![u]\!][\![v]\!]\, ds + \sum_{F \in \mathcal{F}} \frac{\delta_1}{h_i} \int_F [\![\nabla u]\!][\![\nabla v]\!]\, ds \\
& - \sum_{F_i \in \mathcal{F}_D} (\delta_0 v - \mathbf{n}_i \cdot \nabla \Delta v) g_0\, ds - \sum_{F_i \in \mathcal{F}_D} (\delta_1 \mathbf{n}_i \cdot \nabla v + \Delta v) g_1\, ds, \quad \forall v \in H^4(\Omega, \mathcal{T}_h).
\end{aligned}$$

Using Lemma 1, we complete the proof. □

## 4. Analysis of the dG-IgA Scheme

Let $V_h \subset H^4(\Omega, \mathcal{T}_h)$ be the B-spline space defined as

$$V_h := \{v \in L_2(\Omega) : v|_{\Omega_i} \in \mathbb{V}_{h,i}, i = 1, \ldots, N\}, \tag{22}$$

where $\mathbb{V}_{h,i}$ is the B-spline space corresponding to patch $\Omega_i, i = 1, \ldots, N$ for B-splines of degree $p_i \geq 3$. The discrete dG-IgA scheme then reads as: find $u_h \in V_h$ such that

$$a_h(u_h, v_h) = \ell(v_h), \quad \forall v_h \in V_h. \tag{23}$$

An immediate consequence of Theorem 1 is the Galerkin orthogonality property

$$a_h(u - u_h, v_h) = 0, \quad \forall v_h \in V_h. \tag{24}$$

Next, for $v \in H^2(\Omega, \mathcal{T}_h)$, we show that the bilinear form $a_h(\cdot, \cdot)$, is coercive with respect to the following norm

$$\|v\|_h^2 := \sum_{i=1}^N \|\Delta v\|_{L_2(\Omega_i)}^2 + \sum_{F \in \mathcal{F}} \frac{\delta_0}{h_i^3} \|[\![v]\!]\|_{L_2(F)}^2 + \sum_{F \in \mathcal{F}} \frac{\delta_1}{h_i} \|[\![\nabla v]\!]\|_{L_2(F)}^2. \tag{25}$$

**Remark 2.** The discrete norm (25) is indeed a norm. For some function $v \in H^2(\Omega, \mathcal{T}_h)$, if $\|v\|_h = 0$, then

$$\begin{aligned} \Delta v &= 0 \quad \text{in} \quad \Omega_i \quad \text{for all} \quad i = 1, \ldots, N \\ [\![v]\!] &= 0, [\![\nabla v]\!] = 0 \quad on \quad F \in \mathcal{F}. \end{aligned} \tag{26}$$

This means $v = 0$ in the whole domain $\overline{\Omega}$ since it is a weak solution to (26). The homogeneity and triangle inequality axioms are easily verifiable.

**Lemma 2.** *For arbitrary positive numbers $\delta_0$ and $\delta_1$, the discrete bilinear form $a_h(\cdot, \cdot) : V_h \times V_h \to \mathbb{R}$ defined in (23) is $V_h$−coercive with respect to the norm $\|\cdot\|_h$, i.e.*

$$a_h(v_h, v_h) = \|v_h\|_h^2, \quad \forall v_h \in V_h. \tag{27}$$

PROOF. The proof follows by choosing $u_h = v_h$ in the bilinear form (23). □

Lemma 2 yields that the solution of the discrete variational problem (23) is unique. Since the discrete variational problem is in the finite dimensional space $V_h$, the uniqueness yields the existence of the solution $u_h \in V_h$ of (23).

**Remark 3.** The bilinear form $a_h(\cdot, \cdot)$ is meaningful for B-spline or NURBS degree $p_i \geq 3$. We note that for B-splines or NURBS degree $p_i \geq 2$, we obtain the so called continuous/discontinuous Galerkin method for the biharmonic problem, see e.g. [12] as follows; find $u_h \in V_h$ such that

$$\widetilde{a}_h(u_h, v_h) = \widetilde{\ell}_h(v_h), \quad \forall v_h \in V_h, \tag{28}$$

where the bilinear form is given as

$$\widetilde{a}_h(u_h, v_h) = \sum_{i=1}^N \int_{\Omega_i} \Delta u_h \Delta v_h \, dx - \sum_{F \in \mathcal{F}} \int_F [\![\nabla v_h]\!]\{\Delta u_h\} \, ds + \sum_{F \in \mathcal{F}} \int_F [\![\nabla u_h]\!]\{\Delta v_h\} \, ds$$

$$+ \sum_{F \in \mathcal{F}} \frac{\delta_0}{h_i^3} \int_F [\![u_h]\!][\![v_h]\!] \, ds + \sum_{F \in \mathcal{F}} \frac{\delta_1}{h_i} \int_F [\![\nabla u_h]\!][\![\nabla v_h]\!] \, ds, \tag{29}$$

and linear form

$$\widetilde{\ell}_h(v_h) = \sum_{i=1}^{N} \int_{\Omega_i} f v_h \, dx + \sum_{F_i \in \mathcal{F}_D} \delta_0 v_h g_0 \, ds + \sum_{F_i \in \mathcal{F}_D} (\delta_1 \mathbf{n}_i \cdot \nabla v_h + \Delta v_h) g_1 \, ds. \tag{30}$$

The bilinear form (29) is coercive with respect the norm $\|\cdot\|_h$ i.e.

$$\widetilde{a}_h(v_h, v_h) = \|v_h\|_h^2, \quad \forall v_h \in V_h. \tag{31}$$

The discrete stability result (31) yields the uniqueness and existence of the discrete solution $u_h$.

An often necessary tool for the analysis of the dG-IgA scheme are the multi-patch inverse and trace inequalities given by the following lemmata.

**Lemma 3.** *Let $K_i \in \mathcal{K}_{h,i}, i = 1, \ldots, N$. Then the inverse inequalities,*

$$\|\nabla v\|_{L_2(\Omega_i)} \leq C_{inv,1,u} h_i^{-1} \|v\|_{L_2(\Omega_i)}, \tag{32}$$

$$\|v\|_{L_2(\partial \Omega_i)} \leq C_{inv,0,u} h_i^{-1/2} \|v\|_{L_2(\Omega_i)}, \tag{33}$$

*hold for all $v \in V_h$, where $C_{inv,1,u}$ and $C_{inv,0,u}$ are positive constants, which are independent of $h_i$ and $\Omega_i$.*

**Lemma 4.** *Let $K \in \mathcal{K}_{h,i}, i = 1, \ldots, N$ and $\widehat{K} = \Phi_i^{-1}(K)$. Then the scaled trace inequality*

$$\|v\|_{L_2(\partial \Omega_i)} \leq C_{tr,u} h_i^{-1/2} \left( \|v\|_{L_2(\Omega_i)} + h_i |v|_{H^1(\Omega_i)} \right), \tag{34}$$

*holds for all $v \in H^1(\Omega_i)$, where $h_i$ denotes the global mesh size of patch $\Omega_i$ in the physical domain, and $C_{tr,u}$ is a positive constant that only depends on the quasi-uniformity and shape regularity of the mapping $\Phi_i$.*

The proof of the above lemmata (i.e. Lemma 3 and Lemma 4) follows by using the patch-wise inverse and trace inequalities, see e.g., [13] and the quasi-uniformity assumption (8), see [21, chapter 2].

Next, we show that the bilinear form $a_h(\cdot, \cdot)$ $V_{h,*} \times V_h$, is uniformly bounded with $V_{h,*} := V_D \cap H^{\mathbf{s}}(\Omega, \mathcal{T}_h) + V_h$ with $\mathbf{s} \geq 4$ and equipped with the norm

$$\|v\|_{h,*}^2 = \|v\|_h^2 + \sum_{F \in \mathcal{F}} \frac{h_i^3}{\delta_0} \|\{\nabla \Delta v\}\|_{L_2(F)}^2 + \sum_{F \in \mathcal{F}} \frac{h_i}{\delta_1} \|\{\Delta v\}\|_{L_2(F)}^2. \tag{35}$$

The norm (25) is also a norm on $H^4(\Omega, \mathcal{T}_h)$ since $H^4(\Omega, \mathcal{T}_h) \subset H^2(\Omega, \mathcal{T}_h)$.

**Lemma 5.** *Let $a_h(\cdot, \cdot) : V_{h,*} \times V_h$ be the bilinear form defined in (20), then there exists a positive constant $\mu_b$, such that*

$$|a_h(u, v_h)| \leq \mu_b \|u\|_{h,*} \|v_h\|_h, \quad \forall u \in V_{h,*}, v_h \in V_h. \tag{36}$$

PROOF. We estimate the terms in the bilinear form (20) by using the Cauchy-Schwarz inequality.

$$\left|\sum_{i=1}^{N}\int_{\Omega_i}\Delta u\Delta v_h\,dx\right| \leq \left(\sum_{i=1}^{N}\|\Delta u\|^2_{L_2(\Omega_i)}\right)^{1/2}\left(\sum_{i=1}^{N}\|\Delta v_h\|^2_{L_2(\Omega_i)}\right)^{1/2},$$

$$\left|\sum_{F\in\mathcal{F}}\int_F \frac{\delta_0}{h_i^3}[\![u]\!][\![v_h]\!]\,ds\right| \leq \left(\sum_{F\in\mathcal{F}}\frac{\delta_0}{h_i^3}\|[\![u]\!]\|^2_{L_2(F)}\right)^{1/2}\left(\sum_{F\in\mathcal{F}}\frac{\delta_0}{h_i^3}\|[\![v_h]\!]\|^2_{L_2(F)}\right)^{1/2},$$

$$\left|\sum_{F\in\mathcal{F}}\int_F \frac{\delta_1}{h_i}[\![\nabla u]\!][\![\nabla v_h]\!]\,ds\right| \leq \left(\sum_{F\in\mathcal{F}}\frac{\delta_1}{h_i}\|[\![\nabla u]\!]\|^2_{L_2(F)}\right)^{1/2}\left(\sum_{F\in\mathcal{F}}\frac{\delta_1}{h_i}\|[\![\nabla v_h]\!]\|^2_{L_2(F)}\right)^{1/2},$$

$$\left|\sum_{F\in\mathcal{F}}\int_F \{\Delta u\}[\![\nabla v_h]\!]\,ds\right| \leq \left(\sum_{F\in\mathcal{F}}\frac{h_i}{\delta_1}\|\{\Delta u\}\|^2_{L_2(F)}\right)^{1/2}\left(\sum_{F\in\mathcal{F}}\frac{\delta_1}{h_i}\|[\![\nabla v_h]\!]\|^2_{L_2(F)}\right)^{1/2}$$

$$\left|\sum_{F\in\mathcal{F}}\int_F \{\nabla\Delta u\}[\![v_h]\!]\,ds\right| \leq \left(\sum_{F\in\mathcal{F}}\frac{h_i^3}{\delta_0}\|\{\nabla\Delta u\}\|^2_{L_2(F)}\right)^{1/2}\left(\sum_{F\in\mathcal{F}}\frac{\delta_0}{h_i^3}\|[\![v_h]\!]\|^2_{L_2(F)}\right)^{1/2}.$$

Then, concerning the third term, we apply the Cauchy-Schwarz inequality and inverse inequality (33) for $v_h \in V_h$ yielding

$$\left|\sum_{F\in\mathcal{F}}\int_F \{\Delta v_h\}[\![\nabla u]\!]\,ds\right| \leq \left(\sum_{i=1}^{N}\frac{C^2_{inv,0,u}}{\delta_1}\|\Delta v_h\|^2_{L_2(\Omega_i)}\right)^{\frac{1}{2}}\left(\sum_{F\in\mathcal{F}}\frac{\delta_1}{h_i}\|[\![\nabla u]\!]\|^2_{L_2(F)}\right)^{\frac{1}{2}}.$$

Again by using the Cauchy-Schwarz and inverse inequalities (32) and (33) for $v_h \in V_h$, the fifth term in the bilinear form (20) yields

$$\left|\sum_{F\in\mathcal{F}}\int_F \{\nabla\Delta v_h\}[\![u]\!]\,ds\right| \leq \left(\sum_{i=1}^{N}\frac{C^2_{inv,0,1}}{\delta_0}\|\Delta v_h\|^2_{L_2(\Omega_i)}\right)^{\frac{1}{2}}\left(\sum_{F\in\mathcal{F}}\int_F \frac{\delta_0}{h_i^3}\|[\![u]\!]\|^2_{L_2(F)}\right)^{\frac{1}{2}},$$
(37)

where $C^2_{inv,0,1} = C^2_{inv,0,u}C^2_{inv,1,u}$. Finally, by using Cauchy-Schwarz's inequality, we have

$$|a_h(u,v_h)| \leq \Bigg[\sum_{i=1}^{N}\|\Delta u\|^2_{L_2(\Omega_i)} + 2\sum_{F\in\mathcal{F}}\frac{\delta_0}{h_i^3}\|[\![u]\!]\|^2_{L_2(F)} + 2\sum_{F\in\mathcal{F}}\frac{\delta_1}{h_i}\|[\![\nabla u]\!]\|^2_{L_2(F)}$$
$$+ \sum_{F\in\mathcal{F}}\frac{h_i^3}{\delta_0}\|\{\nabla\Delta u\}\|^2_{L_2(F)} + \sum_{F\in\mathcal{F}}\frac{h_i}{\delta_1}\|\{\Delta u\}\|^2_{L_2(F)},\Bigg]^{1/2}$$
$$\times \Bigg[(1+C^2_{inv,0,1}/\delta_0 + C^2_{inv,0,u}/\delta_1)\sum_{i=1}^{N}\|\Delta v_h\|^2_{L_2(\Omega_i)}$$
$$+ 2\sum_{F\in\mathcal{F}}\frac{\delta_0}{h_i^3}\|[\![v_h]\!]\|^2_{L_2(F)} + 2\sum_{F\in\mathcal{F}}\frac{\delta_1}{h_i}\|[\![\nabla v_h]\!]\|^2_{L_2(F)}\Bigg]^{1/2}$$

$$\leq \mu_b \|u\|_{h,*} \|v_h\|_h,$$

with $\mu_b = 2\sqrt{\max\{1, (1 + C_{inv,0,1}^2/\delta_0 + C_{inv,0,u}^2/\delta_1)\}}.$ □

**Lemma 6.** *The norms $\|\cdot\|_h$ and $\|\cdot\|_{h,*}$ are uniformly equivalent on the discrete space $V_h$ such that*
$$\|v_h\|_{h,*} \leq C_* \|v_h\|_h, \quad \forall v_h \in V_h, \tag{38}$$
*where $C_*$ is mesh independent.*

PROOF. By choosing $v_h \in V_h$ and using (35), we have

$$\|v_h\|_{h,*}^2 = \|v_h\|_h^2 + \sum_{F \in \mathcal{F}} \frac{h_i^3}{\delta_0} \|\{\nabla \Delta v_h\}\|_{L_2(F)}^2 + \sum_{F \in \mathcal{F}} \frac{h_i}{\delta_1} \|\{\Delta v_h\}\|_{L_2(F)}^2. \tag{39}$$

Applying the inverse inequality (32) to the second term and again (32) together with (33) to the last term yields

$$\|v_h\|_{h,*}^2 \leq \left(1 + C_{inv,0,1}^2/\delta_0 + C_{inv,0,u}^2/\delta_1\right) \|v_h\|_h^2,$$

where $C_{inv,0,1}^2 = C_{inv,0,u}^2 C_{inv,1,u}^2$. □

Using Lemma 6, the boundedness of the bilinear form $a_h(\cdot,\cdot)$ i.e. Lemma 36, yields

$$|a_h(u_h, v_h)| \leq \tilde{\mu}_b \|u_h\|_h \|v_h\|_h, \quad \forall u_h, v_h \in V_h, \tag{40}$$

where $\tilde{\mu}_b = \mu_b \left(1 + C_{inv,0,1}^2/\delta_0 + C_{inv,0,u}^2/\delta_1\right)^{1/2}.$

## 5. Error Analysis of the dG-IgA Discretization

In this section, we present the approximation estimates required to obtain *a priori* error estimates. For patch $\Omega_i, i = 1, \ldots, N$, let $\Pi_h : L_2(\Omega_i) \to V_{h,i}$ denote a quasi-interpolant that yields optimal approximation results. Of course, such an interpolant is known to exist and has been well studied and presented in [3, 4] as follows

**Lemma 7.** *Let $l_i$ and $s_i$ be integers with $0 \leq l_i \leq s_i \leq p_i + 1$ and $K \in \mathcal{K}_{h,i}$. Then there exist an interpolant $\Pi_h v \in V_{h,i}$ for all $v \in H^{s_i}(\Omega_i)$ and a constant $C_s > 0$ such that the following inequality holds*

$$\sum_{K \in \mathcal{K}_{h,i}} |v - \Pi_h v|_{H^l(K)}^2 \leq C_s h_i^{2(s_i - l_i)} \|v\|_{H^{s_i}(\Omega_i)}^2, \tag{41}$$

*where $h_i$ is the mesh size in the physical domain, and $p$ denotes the underlying polynomial degree of the B-spline or NURBS.*

For patch $\Omega_i, i = 1, \ldots, N$, the local estimate (41) yields a global estimate if the multiplicity of the inner knots is not larger than $p_i + 1 - l_i$ and $\Pi_h v \in V_{h,i} \cap H^{l_i}(\Omega_i)$.

**Proposition 1.** Let us assume that the multiplicity of the inner knots is not larger than $p_i+1-l_i$. Given the integers $l_i$ and $s_i$ such that $0 \leq l_i \leq s_i \leq p_i+1$, there exist a positive constant $C_s$ such that for a function $v \in H^{s_i}(\Omega_i)$

$$|v - \Pi_h v|_{H^{l_i}(\Omega_i)} \leq C_s h_i^{(s_i - l_i)} \|v\|_{H^{s_i}(\Omega_i)}, \tag{42}$$

where $h_i$ denotes the maximum mesh-size parameter in the physical domain and the generic constant $C_s$ only depends on $l_i, s_i$ and $p_i$, the shape regularity of the physical domain $\Omega_i$ described by the mapping $\Phi_i$ and, in particular, $\nabla \Phi_i$.

PROOF. See [26, Proposition 3.2].

For the current analysis, we consider that the quasi-interpolant is the same for each patch, i.e. $\Pi_h : V_D \cap H^{\mathbf{s}}(\Omega, \mathcal{T}_h) \to V_h$ with $\mathbf{s} \geq 4$.

**Lemma 8.** Let $v \in V_D \cap H^{\mathbf{s}}(\Omega, \mathcal{T}_h)$ with a positive integer $\mathbf{s} \geq 4$ and let $F \in \mathcal{F}_I \cup \mathcal{F}_D$ be the facets. By assuming quasi-uniform meshes, then there exists a quasi-interpolant $\Pi_h$ such that $\Pi_h v \in V_h$ and the following estimates hold

$$\|\nabla^q(v - \Pi_h v)\|_{L_2(\partial \Omega_i)}^2 \leq C_0 h_i^{2(r_i - q) - 1} \|v\|_{H^{r_i}(\Omega_i)}^2, \tag{43}$$

$$\sum_{F \in \mathcal{F}} \frac{\delta_0}{h_i^3} \|[\![v - \Pi_h v]\!]\|_{L_2(F)}^2 \leq C_1 \sum_{i=1}^{N} h_i^{2(r_i - 2)} \|v\|_{H^{r_i}(\Omega_i)}^2, \tag{44}$$

$$\sum_{F \in \mathcal{F}} \frac{\delta_1}{h_i} \|[\![\nabla(v - \Pi_h v)]\!]\|_{L_2(F)}^2 \leq C_2 \sum_{i=1}^{N} h_i^{2(r_i - 2)} \|v\|_{H^{r_i}(\Omega_i)}^2, \tag{45}$$

$$\sum_{F \in \mathcal{F}} \frac{h_i}{\delta_1} \|\{\Delta(v - \Pi_h v)\}\|_{L_2(F)}^2 \leq C_3 \sum_{i=1}^{N} h_i^{2(r_i - 2)} \|v\|_{H^{r_i}(\Omega_i)}^2, \tag{46}$$

$$\sum_{F \in \mathcal{F}} \frac{h_i^3}{\delta_0} \|\{\nabla \Delta(v - \Pi_h v)\}\|_{L_2(F)}^2 \leq C_4 \sum_{i=1}^{N} h_i^{2(r_i - 2)} \|v\|_{H^{r_i}(\Omega_i)}^2, \tag{47}$$

where $q$ is a positive integer, $r_i = \min\{s_i, p_i + 1\}$ and the generic constants $C_0, C_1, C_2, C_3$ and $C_4$ are independent of the mesh size.

PROOF. The proof of the first inequality (43) follows by using the trace inequality (34) and the approximation estimate (42) as follows

$$\|\nabla^q(v - \Pi_h v)\|_{L_2(\partial \Omega_i)}^2 \leq C_{tr,u}^2 h_i^{-1} \left( |v - \Pi_h v|_{H^q(\Omega_i)}^2 + h_i^2 |v - \Pi_h v|_{H^{q+1}(\Omega_i)}^2 \right)$$

$$\leq C_0 h_i^{2(r_i - q) - 1} \|v\|_{H^{r_i}(\Omega_i)}^2, \quad \text{with } C_0 = 2C_s C_{tr,u}^2.$$

Next, by using definition (11) and setting $q = 0$ in (43), we have

$$\sum_{F \in \mathcal{F}} \frac{\delta_0}{h_i^3} \|[v - \Pi_h v]\|^2_{L_2(F)} = \sum_{F \in \mathcal{F}} \frac{\delta_0}{h_i^3} \|(v - \Pi_{h,i} v) \mathbf{n}_i + (v - \Pi_{h,j} v) \mathbf{n}_j\|^2_{L_2(F)}$$

$$\leq \sum_{i=1}^{N} 2 \frac{\delta_0}{h_i^3} \|v - \Pi_{h,i} v\|^2_{L_2(\partial \Omega_i)}$$

$$\leq \sum_{i=1}^{N} 2 \frac{\delta_0}{h_i^3} C_0 h_i^{2r_i - 1} \|v\|^2_{H^{r_i}(\Omega_i)} = C_1 \sum_{i=1}^{N} h_i^{2(r_i - 2)} \|v\|^2_{H^{r_i}(\Omega_i)}.$$
(48)

The other estimates (45) to (47) follow by using Definition (11) and the estimate (43) with $q = 1, 2, 3$ and $q = 4$. $\square$

To derive *a priori* error estimates, we will need the following estimates.

**Lemma 9.** *Let $v \in V_D \cap H^s(\Omega, \mathcal{T}_h)$ for $\mathbf{s} := (s_1, s_2, \ldots, s_N) \geq 4$ and $\Pi_h v \in V_h$ be a projection. Then, for $\mathbf{p} \geq 3$ the following bound holds*

$$\|v - \Pi_h v\|^2_h \leq C_5 \sum_{i=1}^{N} h_i^{2(r_i - 2)} \|u\|^2_{H^{r_i}(\Omega_i)},$$
(49)

$$\|v - \Pi_h v\|^2_{h,*} \leq C_6 \sum_{i=1}^{N} h_i^{2(r_i - 2)} \|u\|^2_{H^{r_i}(\Omega_i)},$$
(50)

*where $r_i := \min\{s_i, p_i + 1\}$, $p_i$ is the degree of the B-spline and the constants $C_5$ and $C_6$ are independent of mesh sizes $h_i$.*

PROOF. By using the definition of the norm (25) and Lemma 8, we obtain

$$\|v - \Pi_h v\|^2_h$$

$$= \sum_{i=1}^{N} \|\Delta(v - \Pi_h v)\|^2_{L_2(\Omega_i)} + \sum_{F \in \mathcal{F}} \frac{\delta_0}{h_i^3} \|[v - \Pi_h v]\|^2_{L_2(F)} + \sum_{F \in \mathcal{F}} \frac{\delta_1}{h_i} \|[\nabla(v - \Pi_h v)]\|^2_{L_2(F)}$$

$$\leq C_s \sum_{i=1}^{N} h_i^{2(r_i - 2)} \|v\|^2_{H^{r_i}(\Omega_i)} + C_1 \sum_{i=1}^{N} h_i^{2(r_i - 2)} \|v\|^2_{H^{r_i}(\Omega_i)} + C_2 \sum_{i=1}^{N} h_i^{2(r_i - 2)} \|v\|^2_{H^{r_i}(\Omega_i)}$$

$$= C_5 \sum_{i=1}^{N} h_i^{2(r_i - 2)} \|v\|^2_{H^{r_i}(\Omega_i)},$$

where $C_5 = (C_s + C_1 + C_2)$. Next, by using the definition of the norm (35) together with Lemma 8, we have

$$\|v - \Pi_h v\|^2_{h,*}$$

$$= \|v - \Pi_h v\|^2_h + \sum_{F \in \mathcal{F}} \frac{h_i^3}{\delta_0} \|\{\nabla \Delta(v - \Pi_h v)\}\|^2_{L_2(F)} + \sum_{F \in \mathcal{F}} \frac{h_i}{\delta_1} \|\{\Delta(v - \Pi_h v)\}\|^2_{L_2(F)}$$

$$\leq C_5 \sum_{i=1}^{N} h_i^{2(r_i-2)} \|v\|_{H^{r_i}(\Omega_i)}^2 + C_3 \sum_{i=1}^{N} h_i^{2(r_i-2)} \|v\|_{H^{r_i}(\Omega_i)}^2 + C_4 \sum_{i=1}^{N} h_i^{2(r_i-2)} \|v\|_{H^{r_i}(\Omega_i)}^2$$

$$= C_6 \sum_{i=1}^{N} h_i^{2(r_i-2)} \|v\|_{H^{r_i}(\Omega_i)}^2,$$

where $C_6 = (C_5 + C_3 + C_4)$. □

Finally, we present the main result for the paper.

**Theorem 2.** *Let $u \in V_D \cap H^{\mathbf{s}}(\Omega, \mathcal{T}_h)$ with $\mathbf{s} = \{s_i, i = 1, \ldots, N\}, s_i \geq 4$ be the solution of (19) for non-negative real numbers $\delta_0$ and $\delta_1$, and let $u_h \in V_h$ be the solution of (23). Then there exists $C > 0$ independent of $h_i$ and $N$ such that the following bound holds:*

$$\|u - u_h\|_h \leq C \sum_{i=1}^{N} h_i^{r_i-2} \|u\|_{H^{r_i}(\Omega_i)} \quad (51)$$

*where $r_i = \min\{s_i, p_i + 1\}$.*

PROOF. Using the coercivity result of Lemma 2, the Galerkin orthogonality (24) and the boundedness of Lemma 5, we can derive the following estimates

$$\|\Pi_h u - u_h\|_h^2 = a_h(\Pi_h u - u_h, \Pi_h u - u_h) = a_h(\Pi_h u - u, \Pi_h u - u_h)$$
$$\leq \mu_b \|\Pi_h u - u\|_{h,*} \|\Pi_h u - u_h\|_h.$$

Thus, we have

$$\|\Pi_h u - u_h\|_h \leq \mu_b \|\Pi_h u - u\|_{h,*}. \quad (52)$$

Using (52) and the estimates from Lemma 9, we obtain by using triangle inequality the following estimate

$$\|u - u_h\|_h \leq \|u - \Pi_h u\|_h + \mu_b \|\Pi_h u - u\|_{h,*}$$
$$\leq \left(C_5^{1/2} + \mu_b C_6^{1/2}\right) \sum_{i=1}^{N} h_i^{r_i-2} \|u\|_{H^{r_i}(\Omega_i)}, \quad (53)$$

with $C = (C_5^{1/2} + \mu_b C_6^{1/2})$. □

## 6. Numerical Results

We present two examples of the model problem (1). The numerical experiments are performed in G+Smo[15]. We choose the penalty parameters $\delta_0 = \delta_1 = (p+1)(p+2)/2$ where $p$ is the B-spline or NURBS degree, see Remark 1. We compute the rate of convergence on a successive mesh refinement by using the formula $rate := \log_2\left(e_{i+1}/e_i\right)$, where $e_{i+1} = \|u - u_{h,i+1}\|_h$ and $e_i = \|u - u_{h,i}\|_h$. We consider a unit square $(0,1)^2$ consisting of four patches, see Figure 3 (left) and a quarter annulus consisting of two patches, see Figure 5 (left). The underlying meshes of patches are matching. The resulting linear system from the discrete dG-IgA scheme (23) has been solved using SuperLU solver, see [20].

*6.1. Example I*

We consider the homogeneous Dirichlet problem for the biharmonic equation with $f(x,y)$ chosen such that the analytical solution is given by $u(x,y) = sin^2(\pi x) sin^2(\pi y)$. We consider a unit square domain $\Omega = (0,1)^2$ consisting of four patches, see Figure 3 (left). A similar example has been presented using $hp$ dG-FEM in [25]. The corresponding contours of the solution can be seen in Figure 3 (right). In Figure 4, we observe optimal convergence rate i.e. $\mathcal{O}(h^{p-1})$ in the discrete norm $\|\cdot\|_h$ for $2 \leq p \leq 6$. Although we proved the optimal convergence rate for B-spline and NURBS degree $p \geq 3$, we also observe an optimal convergence rates for $p \geq 2$. This is because the discrete norm $\|\cdot\|_h$ still yields the optimal estimate for the continuous/discontinuous Galerkin scheme, see Remark 3.

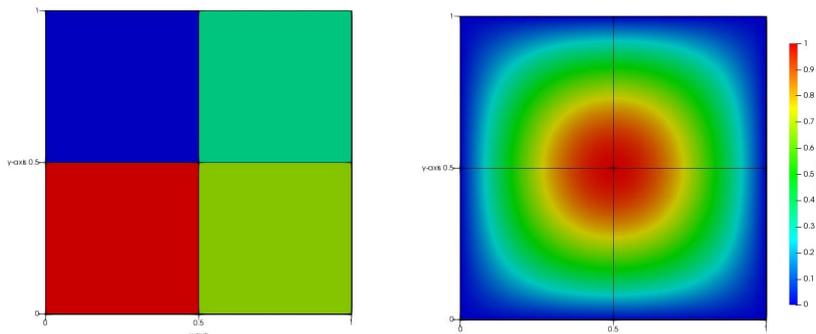

Figure 3: Example of unit square consisting of four patches (left) and solution contours (right).

*6.2. Example II*

We consider again the homogeneous Dirichlet problem for the model problem (1). We choose $f(x,y)$ such that the exact solution is given by $u(x,y) = x^2 y^2 (x^2 + y^2 - 1)^2 (x^2 + y^2 - 4)^2$. We consider a quarter annulus with inner radius $r_{in} = 1$ and outer radius $r_{out} = 4$. The quarter annulus is decomposed into two patches, see e.g. Figure 5 (left). Each patch has Knot vectors $\Xi_1 := \{0,0,1,1\}$ and $\Xi_2 := \{0,0,0,1,1,1\}$. The control points corresponding to patch one (red color) are given by $\mathbf{P}^1_{1,1} = (1.0, 0.0)$, $\mathbf{P}^1_{2,1} = (2.5, 0.0)$, $\mathbf{P}^1_{1,2} = (1.0, 1.0)$, $\mathbf{P}^1_{2,2} = (2.5, 2.5)$, $\mathbf{P}^1_{1,3} = (0.0, 1.0)$ and $\mathbf{P}^1_{2,3} = (0.0, 2, 5)$ and the control points corresponding to patch two (blue color) are given as $\mathbf{P}^2_{1,1} = (2.5, 0.0)$, $\mathbf{P}^2_{2,1} = (4.0, 0.0)$, $\mathbf{P}^2_{1,2} = (2.5, 2.5)$, $\mathbf{P}^2_{2,2} = (4.0, 4.0)$, $\mathbf{P}^2_{1,3} = (0.0, 2.5)$ and $\mathbf{P}^2_{2,3} = (0.0, 4.0)$. The solution contours can be seen in Figure 5 (right).

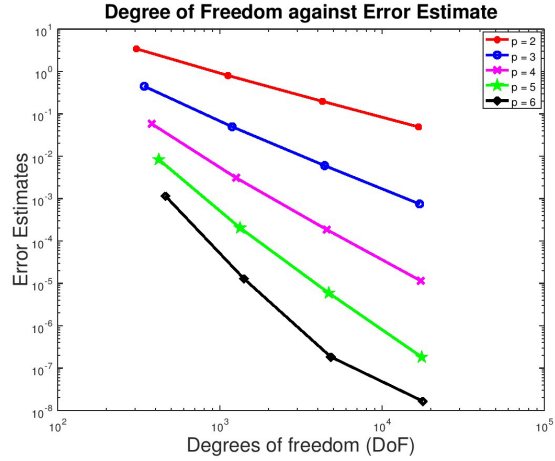

Figure 4: Plot of error estimate in the discrete norm $\|\cdot\|_h$ against the degree of freedom (DoF) for B-splines of degrees $p = 2, 3, 4, 5$ and $6$.

In Figure 6, we observe optimal convergence rate i.e. $\mathcal{O}(h^{p-1})$ in the discrete norm $\|\cdot\|_h$ for $2 \leq p \leq 6$. Although we proved the optimal convergence rate for B-spline and NURBS degree $p \geq 3$, we also observe an optimal convergence rates for $p \geq 2$. This is because the discrete norm $\|\cdot\|_h$ still yields the optimal estimate for the continuous/discontinuous Galerkin scheme, see Remark 3.

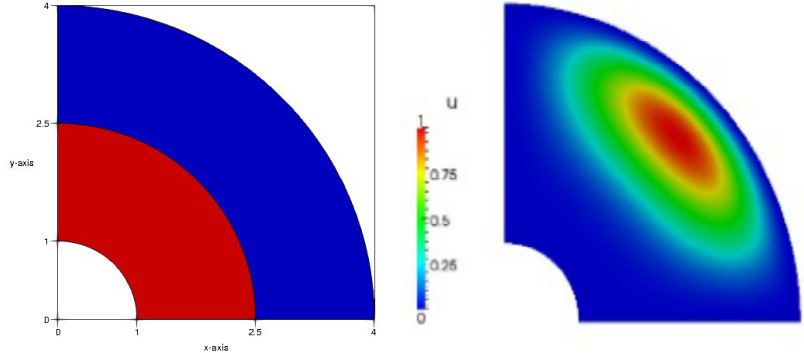

Figure 5: Example of quarter annulus consisting of two patches (left) and solution contours (right).

## 7. Conclusion

In this paper, we presented *a priori* error estimate for the multi-patch discontinuous Galerkin isogeometric analysis (dG-IgA) for the biharmonic prob-lem. We assumed that the solution could be discontinuous on the interior facets and applied interior penalty Galerkin methodology. We showed optimal *a priori* error estimates with respect to a discrete norm $\|\cdot\|_h$ and presented numerical results that confirmed the analysis presented. We have presented the non-symmetric interior penalty Galerkin (NIPG) scheme for the Biharmonic problem. In a forthcoming article, the analysis of the other interior penalty penalty Galerkin schemes for the Biharmonic problem including the symmetric scheme, see e.g. [2] and semi-symmetric schemes will be presented. We will extend the current results to open and closed surfaces following our previous results in [17].

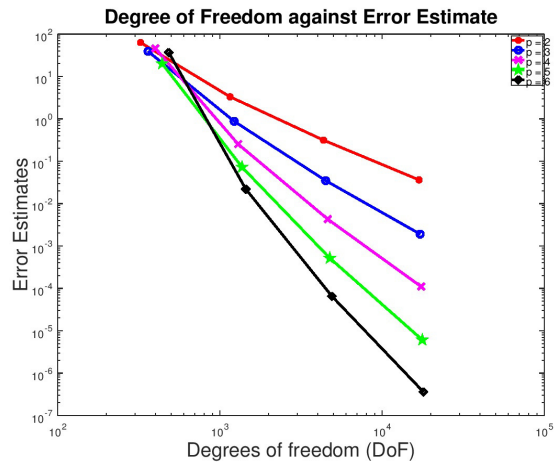

Figure 6: Plot of error estimate in the discrete norm $\|\cdot\|_h$ against the degree of freedom (DoF) for B-splines of degrees $p = 2, 3, 4, 5$ and $6$.

## Acknowledgement

The research is partly supported by the National Institute for Mathematical Sciences (NIMS), Ghana.